\documentclass[11pt]{article}
\usepackage{amsfonts,amssymb,amsmath,a4wide}
\def\Zm{{\mathbb Z}}

\def\Rm{{\mathbb R}}
\def\Tm{{\mathbb T}}
\def\Nm{{\mathbb N}}

\def\lto{\longrightarrow}
\def\lmto{\longmapsto}
\def\eq{\Longleftrightarrow}
\def\leq{\leqslant}
\def\geq{\geqslant}
\newcommand{\mL}{\ensuremath{\mathcal{L}}}
\newcommand{\mT}{\ensuremath{\mathcal{T}}}

\newcommand{\mN}{\ensuremath{\mathcal{N}}}
\newcommand{\mM}{\ensuremath{\mathcal{M}}}
\newcommand{\mS}{\ensuremath{\mathcal{S}}}
\newcommand{\mB}{\ensuremath{\mathcal{G}}}
\newcommand{\mA}{\ensuremath{\mathcal{A}}}

\newcommand{\mH}{\ensuremath{\mathcal{H}}}

\newcommand{\vs}{\vspace{.2cm}}
\makeatletter
\renewcommand{\section}{\@startsection
{section}
{1}
{0mm}
{-1.2\baselineskip}
{\baselineskip}
{\center \scshape}}

\renewcommand{\subsection}{\@startsection
{subsection}
{2}
{0mm}
{-\baselineskip}
{0mm}
{\normalfont \normalsize \bfseries}}
\makeatother

\author{Patrick Bernard}

\title{}

\begin{document}

\begin{center}
\begin{scshape}
\begin{Large}
Connecting orbits of time dependent 
Lagrangian systems
\vspace{1cm}\\
\end{Large}
\begin {large}
Patrick Bernard
\end{large}
\end{scshape}
\vspace{.5cm}\\
October 2001
\end{center}

\begin{itshape}
R\'esum\'e:
On donne une g\'en\'eralisation
\`a la dimension sup\'erieure des r\'esultats 
obtenus par Birkhoff et Mather sur 
l'existence d'orbites errant dans les 
zones d'instabilit\'e des applications de l'anneau
d\'eviant la verticale.
Notre g\'en\'eralisation s'inspire fortement
de celle propos\'ee par Mather dans \cite{Mather}.
Elle pr\'esente cependant l'avantage de contenir effectivement
l'essentiel des
resultats de Birkhoff et Mather sur les diff\'eomorphismes
de l'anneau.
\vspace{.5cm}

Abstract: We generalize to higher dimension 
results of Birkhoff and Mather on the existence 
of orbits wandering in regions of instability of twist maps.
This generalization is strongly
inspired by the one proposed by Mather in 
\cite{Mather}.
However, its advantage is that it 
contains most of the results of Birkhoff and Mather
on twist maps.

\end{itshape}

\vspace{1cm}

A very natural class of problems  in dynamical
systems is  the existence of orbits
connecting prescribed regions of phase space.
There are several important open questions in this
line, like the one posed by Arnold :
Is a generic Hamiltonian system transitive on its 
energy shells?

Birkhoff's theory   of regions of instability of twists maps,
recently extended by Mather using variational methods
and by Le Calvez,
provide very relevant results in that direction.
In short, these works establish the existence,
for a certain class of mappings of the annulus,
of orbits visiting in turn prescribed regions of the annulus
under the hypothesis that these regions are not separated by
a rotational invariant circle.

John Mather has opened the way to a
generalization 
in higher dimension of this
celebrated theory  by proposing what seems to 
be the appropriate setting \textit{i.e.}  time dependent
positive definite Lagrangian systems.
In this setting, he has obtained
the existence of families
of invariant sets generalizing the 
well known Aubry-Mather invariant sets of twist maps.
Then he stated in 1993 a result on the existence of
orbits visiting in turn 
neighborhoods of 
an arbitrary sequence of these invariant sets.
However, the work of Mather
is not a complete achievement since there are no
relevant example in high dimension to which it
can be applied,
and since it is not completely optimal even in the
case of Twist maps.
There are examples where two Aubry-Mather
sets 
of a twist map
are not separated by a rotational 
invariant circle, hence can be connected  by an orbit,
but where this can't be seen by the result of Mather.

In the present paper, we 
state  a new result on the existence
of connecting orbits in higher dimension,
with a full self-contained proof.
This result is very close to the one of Mather,
and the main ideas of the proof are the ones he introduced.
Our result has the advantage that it
is optimal when applied to the twist map case,
but it does not contain the result of Mather,
which we were not able to prove.
\footnote{As it is written in \cite{Mather 2},
the proof contains a gap which I am not able to fill.}

It is still an open question whether these results
may be applied to interesting example in higher dimension
\footnote{
Just before I finished this 
text, John Mather has announced that he had been
able to prove a great result on Arnold diffusion,
so the full achievement of the method may soon be reached.}.
On one hand, it is encouraging that this result
is optimal when restricted to the case
of twist maps,
but on the other hand we will prove that the result 
is useless in the autonomous case.
Additional work will be required 
both to weaken the abstract hypotheses needed  to prove the
existence of connections, 
and to understand when these hypotheses
are satisfied.

\subsection{}\label{hypotheses}
Let $M$ be a smooth, compact, connected manifold, 
$TM\overset{\pi}{\lto}M$ its tangent bundle.
We choose once and for all a Riemaniann metric $g$ on $M$.
It is classical that there is a canonical way to associate to 
it a metric on $TM$.
Let us fix a  $C^2$ Lagrangian function
$
L:TM\times \Rm\lto \Rm.
$
Given any compact interval $I$, we have an action  functional
defined on $C^1(I,M)$ by
$$
A(\gamma)=\int_I L(d\gamma(t),t)dt.
$$
Here and in the following, we note
$d\gamma(t)$ for the curve $d\gamma_t(1):I\lto TM.$
The extremals of $L$ on $I$ are the critical points
of $A$ with fixed endpoints.
We want to study the Lagrangian system associated
with $L$, that is the extremal curves of $L$.
We suppose that $L$ satisfies the following conditions 
introduced by Mather \cite{Mather 1}:\\
\textsc{Periodicity :}
The Lagrangian $L$ is 1-periodic in time \textit{i.e.}
$L(z,t)=L(z,t+1)$ for all $z\in TM$ and all $t\in \Rm$.\\
\textsc{Positive Definiteness :}
For each $x\in M$ and each $t\in \Rm$, the restriction of L to
$T_x M\times t$ is strictly convex with non degenerate Hessian.\\
\textsc{Superlinear Growth :}
For each $t\in \Rm$, 
$$
L(z,t)/\|z\| \lto \infty \,\text{ as } \|z\|\lto \infty.
$$
Under these hypotheses, 
there exists a continuous vector field $E_L$
on $TM\times S$,
the Euler-Lagrange vector-field,
which has the property that a $C^1$ curve
$\gamma$
is an extremal of $L$ if and only if 
the curve $(d\gamma(t),t \text{ mod 1})$ is
an integral curve of $E_L$. Although 
this vector field is only continuous,
it has a flow $\phi_t$ on $TM\times S$
 called the Euler-Lagrange flow. 
We assume :\\
\textsc{Completeness :}
The flow $\phi_t$ is complete \textit{i.e.}
any trajectory $X:I\lto TM\times S$
of the flow can be extended to a trajectory 
$\bar X:\Rm\lto TM\times S.$
\subsection{}
Let $I=[a,b]$ be a compact interval of time.
A curve $\gamma\in  C^1(I,M)$ 
is called a minimizer 
or a minimal curve if it is minimizing the action
among all curves
$\xi \in  C^1(I,M)$
which satisfy
$\gamma(a) =  \xi(a)$
and 
$\gamma(b) =  \xi(b).$
If $J$ is a non compact interval, the curve
 $\gamma \in C^1(J,M)$ is called a minimizer 
if $\gamma \left| _ I \right. $ is minimal for any 
compact interval $I\subset J$.
An orbit $X(t)$ of $\phi_t$ is called minimizing
if the curve $\pi \circ X$ is minimizing,
a point $(z,s)\in TM\times S$ is minimizing
if its orbit $\phi_t((z,s)) $ is minimizing.
Let us call $\tilde \mB$ the set of minimizing points 
of $TM\times S$. We shall see that  $\tilde \mB$
is a nonempty compact subset of  $TM\times S$,
invariant for the Euler-Lagrange flow.

\subsection{}\label{formes}
Let $\eta$ be a  $1$-form of $M\times S$.
We associate to this form a function on
$TM\times \Rm$, still denoted $\eta$,
and defined by
$$
\eta(z, t)=
\langle \eta,(z,t \text{ mod }1,1)
\rangle_{(\pi(z),t\text{ mod }1)},
$$
where $\langle .,.\rangle_{(x,s)}$
is the usual coupling between forms and vectors
of $T_{(x,s)}(M\times S).$
If the form $\eta$ is closed, 
then the Euler-Lagrange vector field 
of $L-\eta$ is the Euler-Lagrange vector field of $L$,
and $L-\eta$ satisfies all the hypotheses of
 \ref{hypotheses}
if $L$ does.
Let us define the mapping 
\begin{align*}
i_s: M & \lto M\times S\\
 x   & \lmto (x,s).
\end{align*}
For any $1$-form $\eta$
on $M\times S$, let us define the
form $\eta_s$ on $M$ by
$$
\eta_s=i_s^*\eta.
$$
If $\eta$ is a closed $1$-form,
we define its class $[\eta]=[\eta_s]\in H^1(M,\Rm)$,
which does not depend on $s$.
Let $\eta$ and $\mu$ be two closed forms such that
$[\eta]=[\mu]$.
It is clear that the minimizing curves of 
$L-\eta$ and $L-\mu$ are the same.
Let us call $\tilde \mB(c)$
the set of minimizing points associated to the Lagrangian
$L-\eta$, where $\eta$ is any closed one-form such that
$[\eta]=c$.
Let us also define,
for each $s\in S$, the set $\tilde \mB_s(c)\subset TM$
of points $z\in TM$ such that $(z,s)\in \tilde \mB(c)$.
We will also call $\mB(c)$ and $\mB_s(c)$
the projections of $\tilde \mB(c)$ and $\tilde \mB_s(c)$
on $M\times S$ and $M$.

\subsection{}\label{ensembleL}
Let $\tilde \omega(c)$ be the 
union of $\omega$-limit points of minimizing 
trajectories $X:[0,\infty)\lto TM\times S.$ 
Let $\tilde \alpha(c)$ be the 
union of $\alpha$-limit points of minimizing 
trajectories $X:(-\infty,0]\lto TM\times S.$ 
In both definitions above, 
minimization is considered with Lagrangians
$L-\eta$, where $\eta$ is any closed one-form on
$M\times \Rm$ satisfying $[\eta]=c.$
We will consider the invariant set 
$$
\tilde \mL(c)=\tilde \omega(c) \cup \tilde \alpha(c).
$$
We will see that $\tilde \mL(c)\subset \tilde \mB(c)$.
In addition, $\tilde \mL$ is contained in the classical 
Aubry set $\tilde \mA(c)$,
 and satisfies the Lipschitz graph property,
see section \ref{globally} for more details.

\subsection{}\label{R'}
We associate to any subset $A$ of $M$ the subspace 
$$
V(A)=
\bigcap\big\{
i_{U*}H_1(U,\Rm) : U \text{ is an open neighborhood of }A \big\}
\subset H_1(M,\Rm),
$$
where $i_{U*}: H_1(U,\Rm) \lto H_1(M,\Rm)$ is the mapping 
induced by the inclusion.
There exists an open neighborhood $U$ of 
$A$ such that 
$V(A)=i_{U*}H_1(U)$.
We can now define, for each $c\in H^1(M,\Rm)$
the following subspace of $ H^1(M,\Rm)$:
$$
R'(c)=
\sum_{t\in S}\left(V\big(\mB_t(c)\big)\right)^{\bot}.
$$
Our improvement compared with \cite{Mather 2} is that 
$R'(c)$ may be bigger than $V\big(\mB_0(c)\big)^{\bot}$, which was 
considered there. 
In fact, the minimizing curves used in Mather's work
satisfy stronger conditions than belonging to $\tilde \mB$,
and their union is a smaller set called the Ma\~ne set 
$\tilde \mN$. As a consequence, our result does not contain the
result stated in \cite{Mather 2}.
However, the proof is only sketched in Mather's paper,
and it is not clear to me how it should be completed.

\subsection {}\label{main}
We say that a continuous curve $c:\Rm\lto  H^1(M,\Rm)$
is admissible if  for each 
$t_0 \in \Rm$, there exists $\delta>0$ such that 
$c(t)-c(t_0)\in R'(c(t_0))$ for all $t\in [t_0-\delta,t_0+\delta]$.
We say $c_0,c_1\in H^1(M,\Rm)$ are C-equivalent
if there exists an admissible  continuous curve $c:[0,1]\lto  H^1(M,\Rm)$
such that $ c(0)=c_0$ and $c(1)=c_1$.
This is precisely the definition of Mather \cite{Mather 2}
except that our $R'(c)$ is different from Mather's one.
We are now in a position to state our 
main result :
\vs \\
\textsc{Theorem : }
Let us fix a C-equivalence class $C$ in $H^1(M,\Rm)$.
Let 
 $(c_i)_{i\in \Zm}$ be a bi-infinite sequence of elements of $C$
and 
$(\epsilon_i)_{i\in \Zm}$ be a bi-infinite sequence of positive 
numbers.
There exist a trajectory $X(t)$ of the 
Euler-Lagrange flow and a bi-infinite increasing sequence 
$t_i$ of times 
such that
$$
d\big(
X(t_i),\tilde \mL(c_i)
\big)
\leq \epsilon_i.
$$
If in addition there exists a class  
$ c_\infty$ such that $c_i= c_\infty$ for large $i$,
or a class $c_{-\infty}$ such that $c_i= c_{-\infty}$ for small $i$,
then the trajectory $X$ 
is $\omega$-asymptotic 
to $\tilde \mL(c_{\infty})$ or 
$\alpha$-asymptotic to 
$\tilde \mL(c_{-\infty})$.
\vs \\
We shall state and prove in section
\ref{connect}
 a slightly refined theorem,
which implies the following corollaries :\\
\textsc{corollary 1 : }
Let $c_0$ and $c_1$ be two C-equivalent classes. There exists 
a trajectory of the Euler Lagrange flow the $\alpha$-limit 
of which lies 
in $\tilde \mL(c_0)$ and the $\omega$-limit of which lies in
$\tilde \mL(c_1)$.\\
\textsc{corollary 2 : }
If there exist two C-equivalent classes  $c_0$ and $c_1$
such that $\tilde \mL(c_0)$ and
$\tilde \mL(c_1)$
are disjoint, then the time one map of the Euler-Lagrange
flow has positive topological entropy.

\subsection{}
Let us insist on the relations between our theorem 
and the theorem of Mather in  \cite{Mather 2}.
The only difference between these two results
lies in the definition of $C$-equivalence,
and more precisely in the definition 
of $R'(c)$.
We replaced 
$$
V\big(\mN_0(c)\big)^{\bot}
$$
as the subspace of allowed directions in \cite{Mather 2}, $\S 12$,
by
$$
R'(c)=
\sum_{t\in S}\left(V\big(\mB_t(c)\big)\right)^{\bot},
$$
where $\mN$ is the set of semi-static curves, see section \ref{globally}.
The bigger the subspace 
of allowed directions is, the stronger the result.
Our result do not contain the result of Mather
because we had to replace the set $\mN$
of semi-static orbits
(see section \ref{globally})
by the larger set $\mB$ of minimizing orbits
in order to fill the proof.
On the other hand our subspace is bigger
in certain cases for example in the twist map case.
An important  consequence is that 
our result is 
optimal in the case $M=S$ while the result of Mather was not.
In this case, two cohomology classes $c$ and $c'$
are $C$-equivalent 
in our sense if and only if
the associated sets $\tilde\mB(c)$
belong to the same region of instability,
that is if they are not separated by an invariant graph.
See section \ref{twist} for the details.
Our result is equivalent to the result of Mather
in the autonomous case, however, as we shall explain in section
\ref{autonomous} it is of no interest in this case.

\subsection{}
In order to apply the theorem,
it is necessary to be able to describe
the $C$-equivalence classes.
This is not an easy task even in the case 
$M=S$.
It requires a good understanding of the set $\mB(c)$
of minimizing curves.
A lot of literature is devoted to the
study of globally minimizing orbits.
We give a review in section \ref{globally}.
We give most of the proofs because most of them
have been written only in the autonomous case.
These results provide a good description of a smaller
set, the Ma\~ne set.
In section \ref{Lax}, we see that 
the difference between the Ma\~ne set and the set $\mB$
is linked with the asymptotic behavior of the so called 
Lax-Oleinik semi-group. 
We exploit this remark to
obtain some results on the shape of the set
$\mB$.
In section \ref{twist},
we apply these results to the case of twist maps,
and obtain that our theorem is optimal in this case.
Unfortunately, there is no hope to apply our
result in the autonomous case,
as is explained in section \ref{autonomous}.

\newpage
\section{Minimization}\label{basic}
It is useful to work in a slightly more general setting.
In this section, 
we will consider a Lagrangian 
$L:TM\times \Rm\lto \Rm$,
not necessarily time-periodic,
satisfying positive definiteness and superlinearity,
but not completeness.

\subsection{}
If the positive definiteness and superlinear growth 
are satisfied, 
there is a continuous flow 
$\psi_t$ on $TM$
such that 
the curve $\gamma$ is a $C^1$ extremal of $L$
if and only if the curve
$X(t)=d\gamma(t)$ is a trajectory of $\psi_t$.
We still call this flow the Euler-Lagrange flow.
This flow is not assumed to be complete in 
the present  section.

\subsection{}
Let $ H\subset H_1 (M,\Rm)$ be the image of the Hurewitcz 
homomorphism, and $K\subset \pi_1(M)$ its kernel.
We shall consider the Abelian covering
$$
\bar   M\overset{p}{\lto}M.
$$
It is the Galois Covering of $M$ which has $K$ as fundamental group.
Its group of deck transformations is canonically isomorphic
to $H$, which is a lattice in $H_1(M,\Rm)$.
In the case $M=\Tm^n$, $\bar M$ is simply the universal cover 
$\Rm^n$.

\subsection{}\label{tonelli}
The variational study of $L$ relies on some standard results
proved in  \cite{Mather 1}. \\
\textsc{Lemma : }
Given a real number $K$ and a compact interval $[a,b]$,
the set
of all absolutely continuous curves 
$\gamma : [a,b] \lto M$
for which $A(\gamma)\leq K$ 
is compact for the topology of uniform convergence.\\
\textsc{Tonelli's theorem : }
Let $[a,b]$ be a compact interval,
and let us fix two points $x_a$ and $x_b$ in $\bar M$.
The action takes a finite minimum  over the set of
absolutely continuous curves 
$\gamma : [a,b] \lto M$ 
which have a lifting $\bar \gamma$ 
satisfying
$ \bar \gamma (a)=x_a$ and $\bar \gamma(b)=x_b$.
If in addition the Euler-Lagrange flow is complete, 
then any curve $\gamma$
realizing this minimum is $C^1$
and $d\gamma(t)$ is a trajectory of
the Euler-Lagrange flow.

Let $I=[a,b]$ be a compact interval of time.
A curve $\gamma\in  C^{ac}(I,M)$ is called a 
$\bar M$-minimizer if one (hence any) of its 
liftings $\bar \gamma $ to the cover $\bar M$ is minimizing
the action among all curves 
$\bar \xi \in  C^{ac}(I,\bar M)$
which satisfy
$\bar \gamma(a) = \bar \xi(a)$
and 
$\bar \gamma(b) = \bar \xi(b).$
A curve $\gamma\in  C^{ac}(I,M)$ 
is called a minimizer 
or a minimal curve if it is minimizing the action
among all curves
$\xi \in  C^{ac}(I,M)$
which satisfy
$\gamma(a) =  \xi(a)$
and 
$\gamma(b) =  \xi(b).$
Minimizers are $\bar M$-minimizers.
A curve  $\gamma \in C^{ac}(\Rm,M)$ is called a minimizer 
if $\gamma \left| _ I \right. $ is minimal for any compact interval $I$.
Let us notice that if the completeness is not assumed,
the absolutely continuous minimizers need not be
$C^1$, an example of this is given in 
\cite{BM}.

\subsection{}\label{existence}
\textsc{Proposition :}
There exist absolutely continuous minimizers 
 $\gamma \in C^{ac}(\Rm,M)$. If the flow is complete,
these minimizers are $C^1$ extremals and  the curves
$d\gamma(t)$
are trajectories  of the Euler-Lagrange flow.

This proposition follows from the following lemmas, 
which are stated in higher generality for later use.

\subsection{}\label{apriori}
\textsc{Lemma :}
Let us fix a positive definite superlinear Lagrangian $L$,
a compact interval of time $[a,b]$  
and a positive constant $C$.
There exists a constant $K$ with the following property:
If $\tilde L$ is a positive definite superlinear Lagrangian 
such that
$$
|\tilde L(z,t)-L(z,t)|\leq C(1+\|z\|)
$$
for all $z\in TM$ and all $t \in [a,b]$,
and if $\gamma:[a,b]\lto M$ is a minimizer of 
$\tilde L$,
then 
$$
\int _a ^b \|d\gamma(t)\|\, dt \,\leq K\;
\text{ and } \;
\int _a^b L(d\gamma(t),t)\,dt \,\leq K.
$$

\textsc{Proof :}
There exists a constant $B$ depending on 
$L$,
$C$ and $[a,b]$ such that
all minimizer $\gamma$ of $\tilde L$
satisfies 
$\tilde A(\gamma)\leq B$,
where $\tilde A$ is the action associated to 
$\tilde L$.
Since $L$ is superlinear, there exists a constant $D$ such that
$$
L(z,t)\geq (C+1)\|z\|-D
$$
for all $z\in TM$ and $t\in [a,b]$.
It follows that $\tilde L\geq \|z\| -C-D,$
and we get the first estimate
$$
\int \|d\gamma\|\leq B+(b-a)(C+D).
$$
We get the second estimate thanks to the inequality 
$$
A(\gamma)
\leq \tilde A(\gamma)+C \int \|d\gamma\|+C(b-a).
$$
This ends the proof of the lemma.\hfill

\subsection{}\label{limit}
\textsc{Lemma :}
Let $L$ be a positive definite superlinear 
Lagrangian,
and let $[a,b]$
be a compact interval of time.
Let $L_n$ be a sequence of positive definite superlinear 
Lagrangians, such that 
$|L_n(z,t)-L(z,t)|\leq \epsilon_n (1+\|z\|)$
for all $z\in TM$ and all $t\in [a,b]$,
where $\epsilon_n$ is a sequence converging to $0$.
If $\gamma_n:[a,b]\lto M$
is a sequence of minimizers of $L_n$ 
converging uniformly to 
$\gamma:[a,b]\lto M$,
then 
$$
A(\gamma)=\lim \, \int _a ^b L_n(d\gamma_n(t),t)\, dt
$$ 
and 
$\gamma$ is a minimizer of $L$ on $]a,b[.$

\textsc{Proof :}
In view of Lemma \ref{apriori},
the sequence $A(\gamma_n)$ is bounded
and $A(\gamma_n)-A_n(\gamma_n)\lto 0.$
By
Lemma \ref{tonelli},
the curve $\gamma$ is  absolutely continuous,  and satisfies 
$$
A(\gamma)\leq \liminf A(\gamma_n)=\liminf A_n(\gamma_n).
$$
In order to prove the lemma,
it is thus sufficient to prove that if $x: [a,b]\lto M$ is an 
absolutely continuous curve
such that $\gamma(t)=x(t)$
in a neighborhood of $a$ and $b$, then 
$A(x)\geq \limsup A_n(\gamma_n).$
Let $x(t)$ be such a curve.
Recall that $x$ is differentiable almost everywhere.
Let us consider an interval $[a',b']\subset [a,b]$
such that $x$ is differentiable at $a'$
and $b'$ and such that $\gamma(a')=x(a')$
and $\gamma(b')=x(b')$.
There exist positive constants $\delta_0$ and $K$ such that,
for all $\delta \in ]0,\delta_0[,$
$$
d\big(   x(a'),x(a'+\delta)  \big)    \leq K\delta\;
\text{ and } \;
d\big(    x(b'-\delta),x(b')   \big)    \leq K\delta.
$$
As a consequence, there exists an integer $N(\delta)$
such that 
$$
d\big(   \gamma_n(a'),x(a'+\delta)  \big)    \leq 2K\delta\;
\text{ and } \;
d\big(    x(b'-\delta),\gamma_n(b')   \big)    \leq 2K\delta
$$
for all $n\geq N(\delta)$.
Now let us consider the geodesic
$\xi :[a',a'+\delta]\lto M$ 
connecting $\gamma_n(a')$ and $x(a'+\delta)$,
and the geodesic
$\zeta:[b'-\delta,b']\lto M$ 
connecting
$x(b'-\delta)$ and $\gamma_n(b')$.
If $\delta\leq \delta_0$ and $n\geq N(\delta)$,
they satisfy
$\|d\xi\|\leq 2K$
and $\| d \zeta \| \leq 2K$,
hence there exists a constant $B$ such that
 $A_n(\xi)\leq B\delta$
and $A_n(\zeta)\leq B\delta$.
Since $\gamma_n$ is minimizing on $[a',b']$,
it follows that
$$
A_n  \big(  x|_{[a'+\delta,b'-\delta])}\big)
    +2B\delta
\geq 
A_n   \big(
\gamma_n|_{[a',b']}
\big).
$$
Taking  the limit, we obtain
$$
A \big(  x|_{[a'+\delta,b'-\delta])}\big)
 +2B\delta
\geq 
 \limsup A   \big(
\gamma_n|_{[a',b']}
\big).
$$
since this holds for all $\delta\leq \delta_0$,
we get that 
$
A \big(
\gamma|_{[a',b']}
\big)
\geq \limsup
A \big(
x_n|_{[a',b']}
\big).
$
At the limit $a'\lto a$, $b'\lto b$,
we obtain that
$A(x)\geq \limsup A(\gamma_n).$\hfill

\subsection{}\label{minlim}
\textsc{Lemma :}
Let $I_n=[a_n,b_n]$ be a nondecreasing sequence of compact intervals
and let $J=\cup _n I_n$.
Let $L_n$ be a sequence of positive definite superlinear
Lagrangians, such that
$$
\big|
L_n(z,t)-L(z,t)
\big|
\leq \epsilon_n(1+\|z\|)
$$
for all $z\in TM$ and all $t\in I_n$,
where $\epsilon_n\lto 0$.
If $\gamma_n:I_n\lto M$ is a sequence 
of minimizers of $L_n$,
then there is an absolutely 
continuous curve $\gamma:J\lto M$
which is minimizing for $L$ on the interior of $J$,
and a subsequence of $\gamma_n$
which converges uniformly on compact sets of $J$ to $\gamma$.
\\
\hspace{.2cm}
\textsc{Proof :}
In view of Lemma \ref{apriori},
the sequence 
$$
k\lmto A\big(
\gamma_k|_{I_n}
\big)
$$
is bounded for each $n$.
It follows from Lemma \ref{tonelli}
that there is a subsequence
of $k\lmto \gamma_k|_{I_n}$
converging uniformly.
By diagonal extraction, we can build a subsequence
of $\gamma_n$ which converges uniformly on compact sets
to an absolutely continuous limit $\gamma:J\lto M$.
By Lemma \ref{limit}, 
this limit is a minimizer of $L$ on the interior of $J$.
\hfill

\subsection{}
We will have in the following to consider one-forms
on $M\times \Rm$ which are neither periodic nor closed.
Let $\mu$ be a  $1$-form of $M\times \Rm$.
We associate to this form a function on
$TM\times \Rm$, still denoted $\mu$,
and defined by
$$
\mu(z, t)=
\langle \mu, (z,t,1)\rangle_{(\pi(z),t)}.
$$
The new Lagrangian $L-\mu$ is positive definite
and superlinear if $L$ is. If $\mu$ is closed,
then the Euler-Lagrange flows of $L$ and $L-\mu$ are 
the same.
Let us define the mapping 
\begin{align*}
i_t: M & \lto M\times \Rm\\
 x   & \lmto (x,t),
\end{align*}
and the form $\mu_t=i_t^*\mu.$
If $\mu $ is closed, we define its homology 
$[\mu]=[\mu_t]\in H^1(M,\Rm)$.
We will often identify a form $\eta$ on $M\times S$
with its periodic pull-back on $M\times \Rm$.
\newpage

\section{Connecting orbits}\label{connect}
In this section, we prove Theorem \ref{main}.
In fact, we will prove a more precise result, 
Theorem \ref{main'}, which clearly implies 
Theorem \ref{main} and the corollaries.
We suppose from now on that $L$ satisfies
all the hypotheses of \ref{hypotheses}.

\subsection{ }
\textsc{Proposition :}
The set $\tilde \mB(c)$ as defined in \ref{formes}
is a non empty compact subset of $TM\times S$.
It is invariant under the Euler-Lagrange flow.
The mapping $c\lmto \tilde \mB(c)$
is upper semi-continuous.\\
\textsc{Proof :} That $\tilde \mB(c)$ is not empty follows from
Proposition \ref{existence}.
The other statements are consequences of the following 
lemma.
\subsection{}\label{compacite}
\textsc{Lemma : }
Let us consider a sequence $c_n\lto c$
of cohomology classes, a sequence $T_n\lto \infty$
of times, and a sequence
$\gamma_n :[-T_n,T_n]\lto M$
of curves minimizing $L-c_n$.
Then there exists a curve $\gamma \in C^1(\Rm,M)$
minimizing $L-c$  and
a subsequence $\gamma_k$ of $\gamma_n$
such that the sequence
$d\gamma_k$
is converging uniformly on compact sets to 
$d\gamma.$

\textsc{Proof :}
This lemma is mainly
a special case of Lemma \ref{limit}.
However, we have to prove that the convergence 
of $\gamma_n$ to $\gamma$ holds in $C^1$
topology.
This is a direct consequence of the theorem of Ascoli
and of the following lemma,  proved in \cite{Mather 1},
on pages 182 and 185.

\subsection{}\label{apriori2}
\textsc{Lemma : }
For all $K\geq 0$, there exists $K'\geq 0$
such that,
if $\gamma :[a,b]\lto M$
is a $\bar M$-minimizer all coverings
$\bar \gamma$ of which satisfy
$$
d(\bar \gamma(b),\bar \gamma(a))\leq K(b-a)
$$
then for each $t\in [a,b]$, 
$$
\|d \gamma(t)\|\leq K'.
$$
\textsc{Corollary : }
Let us consider a compact set $Q\subset H^1(M,\Rm)$.
There exists a constant $K'>0$ such that,
if $b\geq a+1$,
all  curve $\gamma : [a,b]\lto M$ minimizing $L+c$,
with any $c\in Q$
satisfy
$ \|d \gamma(t)\|\leq K'$
for each $t$.

\subsection {}\label{Mset}
The restriction of the Euler-Lagrange flow defines a continuous flow
on the compact set $\tilde \mB(c)$. 
By the Krylov Bogolioubov
theorem, this flow has invariant probability measures.
The Mather set $\tilde \mM(c)$
is the closure of the union of all the supports of these invariant
probability measures.
We have the following lemma, which is a straightforward result of 
topological dynamics:\\
\textsc{Lemma : }
For all positive number $\epsilon$, there exists
a positive number $T$ such that, if 
$X:[0,T]\lto \tilde \mB(c)$ is a trajectory of the Euler-Lagrange flow,
there exists a time $t\in [0,T]$ such that 
$d(X(t),\tilde \mM(c))\leq \epsilon$.

\subsection {}\label{step}
Let $U$ be an open subset of $M\times S$.
We also note $U$ the open subset in $M\times \Rm$
of points $(x,t)$ such that $(x,t\text{ mod }1\in U$.
The one from $\mu$ of $M\times \Rm$ is called a $U$-step form
if there exist a closed form $\bar \mu$
on $M\times S$,
also considered  as a periodic one-form on 
$M\times \Rm$,
such that
the restriction of $\mu$ to $t\leq 0$ is $0$, the restriction
of $\mu$ to $t\geq 1$ is $\bar \mu$,
and such that the restriction 
of $\mu$ to 
 the set 
$U\cup \{t\leq 0\}\cup \{t\geq 1\}$
is closed.

\subsection {}\label{resonance}
We define the subset $R(c)$ of $H^1(M,\Rm)$
as follows :
A class $d$ belongs to $R(c)$ if and only if there exist
an open neighborhood $U$ of  $\mB (c)$
and a $U$-step form $\mu$ such that $[\bar\mu]=d$.
Since  $H^1(M,\Rm)$ is finite dimensional,
there exists an open neighborhood $U$ of  $\mB (c)$
such that, for each $d\in R(c)$,
there exists an $U$-step form satisfying 
$[\bar\mu]=d$. Such a neighborhood $U$ will be called
an adapted neighborhood.
Recalling that $R'(c)$ has been defined in \ref{R'},
we have the inclusion:
$$
R'(c)\subset R(c).
$$
\textsc{Proof :} 
It is enough to prove that for each $t$,
$V\big(\mB_t(c)\big)^{\bot}\subset R(c)$.
Let us fix a time $t\in [0,1]$.
There exist an open neighborhood $\Omega$
of $\mB_{t}(c)$ and a $\delta>0$
such that  $V(\Omega)=V(\mB_{t}(c))$
and such that  $\mB_{s}(c)\subset \Omega$ for all
 $s\in[t-\delta,t+\delta]$.
Given any class $d \in V(\Omega)^{\bot}$,
we take a closed 1-form $\bar \mu$
on $M$ the support of which is disjoint from $\Omega$
and such that $[\bar  \mu]=d$.
We can consider this   one-form on $M$
as a form on $M\times S$. 
Let $f:\Rm\lto \Rm$ be a smooth function
such that $f=0$ on $(-\infty,t-\delta]$ and $f=1$
on $[t+\delta,\infty)$. 
It is not hard to see that the form 
$$
\mu =f(t)\bar  \mu
$$ 
is an  $U$-step form satisfying
$[\bar \mu]=d$,
where $U$ is the open set 
$M\times [0,t-\delta[ \;\cup \;
\Omega \times S \; \cup \; M \times ]t+\delta,1]$.

\subsection{ }
\textsc{Proposition : }\label{prop}
Let us fix  a cohomology class $c$ in $H^1(M,\Rm)$,
and let $U$ be an adapted neighborhood of $\mB(c)$.
There exists a  positive numbers $\delta$
and an integer $T_0$
with the following property :
If  $\eta_0$ is a closed one-form of $M \times \Rm$
satisfying
$[\eta_0]=c$
and 
if $d\in R(c)$ satisfies $|d|\leq \delta$,
then there exists
an $U$-step form $\mu$
satisfying $[\bar \mu]=d$
and such that 
all the minimizers $\gamma :[-T_0,T_0+1]\lto M$ of $L-\mu-\eta_0$
are $C^1$ extrema of $L$.
\vs\\
\textsc{Proof }
The minimizers of $L-\eta_0 - \mu$ do not depend
on the choice of the form $\eta_0$ satisfying
$[\eta_0]=c$. As a consequence, it is enough
to prove the proposition for a fixed
form $\eta_0$.
Since $H^1(M,\Rm)$ is finite dimensional,
it is possible to take a finite dimensional subspace
$E$ of the space of all U-steps  forms on $M\times S$
such that the restriction to $E$ of the linear map
$\mu \lmto [\bar \mu]$
is onto.
We shall prove by contradiction that,
if $\mu\in E$ is  sufficiently small,
there exists a minimizer $\gamma :[-T_0,T_0+1]\lto M$
of  $L-\eta_0 - \mu$  such that $(\gamma(t),t)\in U$
for all $t \in [0,1]$.
Else, there would exist a sequence  $\mu_n$
of elements of $E$ such that $\mu_n\lto 0 $
(this is meaningful in the finite dimensional vector space $E$)
and a sequence 
$\gamma_n:[-T_n,T_n +1]\lto M$, with $T_n\lto\infty$,
of absolutely continuous
curves minimizing $L-\eta_0 - \mu_n$, such that 
$(\gamma_n(t_n),t_n)\not \in U$ for some $t_n\in[0,1]$.
There exists a sequence $\epsilon_n$ of positive numbers
such that $\epsilon_n \lto 0$ and
$$
|\mu_n(z,t)|\leq \epsilon_n\|z\|
$$
for all $(z,t)\in TM\times \Rm$.
By Lemma \ref{minlim}, 
there exist a curve $\gamma\in C^1(\Rm,M)$
minimizing for $L-\eta_0$
and a subsequence of $\gamma_n$
converging uniformly on
compact sets to $\gamma.$
This implies that 
$(\gamma_n(t), t
\text{ mod } 1)\in U$
for all $t\in [0,1]$ when $n$ is large enough,
which is a contradiction.
This ends the proof of the existence of a minimizer
$\gamma:[-T_0,T_0+1]\lto M$
of  $L-\eta_0 - \mu$  such that $(\gamma(t),t)\in U$
for all $t \in [0,1]$.
The form $\eta_0+\mu$ is closed in a neighborhood
of the set $\big\{ (\gamma(t),t)_{t\in \Rm}\big \}\subset M\times \Rm,$
hence $\gamma$ is a $C^1$ extremal of $L$.
\hfill

\subsection {}\label{stairs}
We say that a continuous curve $c:\Rm\lto  H^1(M,\Rm)$
is admissible if  for each 
$t_0 \in \Rm$, there exists $\delta>0$ such that 
$c(t)-c(t_0)\in R(c(t_0))$ for all $t\in [t_0-\delta,t_0+\delta]$.
We say $c_0,c_1\in H^1(M,\Rm)$ are C-equivalent
if there exists an admissible  continuous curve $c:\Rm\lto  H^1(M,\Rm)$
such that $ c(t)=c_0$ for all $t\leq 0$  and 
$c(t)=c_1$ for all $t\geq 1$.
This is precisely the definition of Mather \cite{Mather 2}
and of \ref{main}, 
except that our $R(c)$ is now bigger.
\vspace{.2cm}\\
\textsc{Lemma : }
Let $c_0$ and $c_1$ be two C-equivalent classes.
There exist 
an integer $T(c_0,c_1)$ and a form $\mu$ on $M\times \Rm$ such that :

$\iota.\,$
The restriction of $\mu$ to $\{t\leq 0\}$
is $0$
and the  restriction of $\mu$ to $\{t\geq T(c_0,c_1)\}$
is a closed periodic one form $\bar \mu$
satisfying  $[\bar \mu]=c_1-c_0.$

$\iota \iota .\,$
If $\eta_0$ is a closed periodic one form such that 
$[\eta_0]=c_0$, then any absolutely continuous curve
$\gamma : [0,T(c_0,c_1)]\lto M$ minimizing for
$L-\eta_0-\mu$
is an extremal of $L$. 
\vs\\
\textsc{Proof :}
Let $c(t):\Rm \lto H^1(M,\Rm)$ be an admissible curve
such that $c(t)=c_0$ for all $t\leq 0$ and 
$c(t)=c_1$ for all $t\geq 1$.
Let us fix, for each $t\in [0,1]$, an adapted neighborhood 
$U(t)$ of $\mB(c(t))$,
and let $\delta(t)$ and $T_0(t)$
be the numbers given by applying Proposition \ref{prop}
to $c(t)$ and $U(t)$.
For each $t$, there is a positive number $\delta'(t)$ such that
$
c(s)-c(t)\in R(c(t))$ and $|c(s)-c(t)|\leq  \delta(t)
$
for all $t\in ]t-10\delta'(t),t+10\delta'(t)[$.
There is a finite increasing sequence $(t_i)_{0\leq i\leq N}$ of times
such that the intervals $ ]t_i-\delta'(t_i),t_i+\delta'(t_i)[$
cover $[0,1]$. We require in addition
that $t_0=0$ and $t_N=1$.
To sum up,
we have constructed a finite   increasing sequence 
$(t_i)_{0\leq i\leq N}$
such that
$$
c(t_{i+1})-c(t_i)\in R(c(t_i))
\text{ and } 
|c(t_{i+1})-c(t_i)|\leq \delta(t_i).
$$
Let us  call $\mu_i$ the step form given by
Proposition  \ref{prop}  applied with $d=c(t_{i+1})-c(t_i)$ 
for $0\leq i <N$.
Let us set
 $
T_i=1+ \max \big(T_0(t_i),T_0(t_{i+1}) \big),\;\;0\leq i\leq N-1
$
and $T_{-1}=T_0(t_0)+1$
and define the integers  $(\tau_i)_{-1\leq i\leq N}$ by 
$\tau_0=0$ and $\tau_{i+1}=\tau_i+T_i+1$.
We also consider $\tau_i$ as the  time translation
$
(q,t)\lmto (q,t+\tau_i)
$
on $M\times \Rm$.
Let us define the one form
$$
\mu=\sum _{i=0}^{N-1} (-\tau_i)^* \mu_i.
$$
If $\gamma:\Rm \lto M$ is a minimizer of 
$L-\eta_0-\mu$, then $ \gamma $ is an extremum of $L$.
To check this
let us consider, for  each $1\leq i\leq N-1$, the curve 
$$
\gamma (t+\tau_i): [\tau_{i-1}-\tau_i+1,\tau_{i+1}-\tau_i]
\lto M, 
$$
which is a minimizer of 
$$
L-\eta_0-\sum _{j=0}^{i-1} (\tau_i-\tau_j)^* \bar \mu_j
\;-\; \mu_i,
$$
where $\eta_0+\sum _{j=0}^{i-1} (\tau_i-\tau_j)^*\bar \mu_j$
is a closed form satisfying
$$
\left[
\eta_0+\sum _{j=0}^{i-1} (\tau_i-\tau_j)^*\bar \mu_j
\right]
=c(t_i).
$$
Since $\tau_{i-1}-\tau_i+1=-T_{i-1}\leq -T_0(t_i)$
and since $\tau_{i+1}-\tau_i= T_i+1\geq T_0(t_i)$,
we are in a position to apply Proposition \ref{prop}
and obtain that $\gamma$ is an extremum of $L$
on $[\tau_{i-1}+1, \tau_{i+1}]$
for each $i$ satisfying $1\leq i\leq N-1.$
It follows that $L$ is an extremum of $L$
on $[-T_0,\tau_N]$. 
Since $\eta$ is a closed periodic one-form on
each of the intervals
$(\-\infty,0]$ and $[\tau_N-T_{N-1},\infty)$,
$\gamma$ is clearly an extremum of $L$ on 
these intervals, hence on the whole of $\Rm$.
\hfill

\subsection{}\label{sejour}
\textsc{Lemma : }
For each cohomology class $c$ and each positive number 
$\epsilon$,
there exists a positive number $\mT_{\epsilon}(c)$
with the following property : 
If 
$X:[0,\mT_{\epsilon}(c)]\lto TM\times S$
is a trajectory 
of the Euler-Lagrange flow minimizing $L-c$,
then there
exists a time $t$ in $[0,\mT_{\epsilon}(c)]$
such that
$$
d\big(
X(t),\tilde \mM(c)
\big)\leq \epsilon.
$$
\textsc{Proof : }
Let us fix $\epsilon>0$,
and consider a sequence $X_i:[0,2i]\lto TM\times S$ of
trajectories minimizing $L+c$.
By lemma \ref{compacite},
there exists a minimizing trajectory 
$X\in C(\Rm,TM\times S)$
such that the curves 
$Y_k(t)=X_k(t+k)$
are converging uniformly on compact sets to 
$X(t).$
On the other hand, by Lemma \ref{Mset},
there exists
a time $t$ such that 
$$
d\big(
X(t),\tilde \mM(c)
\big)
\leq \epsilon/2.
$$
It follows that 
$$
d\big(
X_k(t+k),\tilde \mM(c)
\big)
\leq \epsilon
$$
when $k$ is large enough.
\hfill

\subsection{ }\label{main'}
\textsc{Theorem : }
Let us fix a C-equivalence class $C$ in $H^1(M,\Rm)$.
Let 
 $(c_i)_{i\in \Zm}$ be a bi-infinite sequence of elements of $C$
and 
$(\epsilon_i)_{i\in \Zm}$ be a bi-infinite sequence of positive 
numbers.
If $t_i'$ and $t''_i$ are  bi-infinite sequences of real numbers such that
$ t''_i-t'_i\geq \mT_{\epsilon_i}(c_i)$
and $t'_{i+1}-t''_i \geq T(c_i, c_{i+1})$,
then there exist a trajectory $X(t)$ of the 
Euler-Lagrange flow and a bi-infinite sequence 
$t_i\in ]t'_i,t''_i[$
such that
$$
d\big(
X(t_i),\tilde \mM(c_i)
\big)
\leq \epsilon_i.
$$
If in addition there exists a class  
$ c_\infty$ such that $c_i= c_\infty$ for large $i$,
or a class $c_{-\infty}$ such that $c_i= c_{-\infty}$ for small $i$,
then the trajectory $X$ 
is $\omega$-asymptotic 
to $\tilde \mL(c_{\infty})$ or 
$\alpha$-asymptotic to 
$\tilde \mL(c_{-\infty})$.
Recall that the sets $\tilde \mL$ have been defined in \ref{ensembleL}.
\vs\\
\textsc{corollary : }
If there exist two C-equivalent classes  $c_0$ and $c_1$
such that $\tilde \mM(c_0)$ and
$\tilde \mM(c_1)$
are disjoint, then the time one map  of the Euler Lagrange
flow has positive topological entropy.
\vs\\
\textsc{Proof : }
The proof will be quite similar to the proof of Lemma
\ref{stairs}.
Using this lemma,
one can build a 1-form
$\eta$ on $M\times \Rm$
 such that the minimizers of $L-\eta$
are extremals of $L$, and such that,  for each $i$,
the form 
$\eta |_ {[t'_i,t''_i]}$
is closed and periodic and satisfies
$$
\big[
\eta | _{[t'_i,t''_i]}
\big]=c_i.
$$
Let us consider 
a minimizer $\gamma(t)$  of $L-\eta$,
and the as\-sociated trajectory of the Euler-Lagrange flow
$X(t)=(d\gamma(t),t \text{ mod } 1).$
By Lemma \ref{sejour}, there exists a sequence
$t_i\in ]t'_i,t''_i[$ of times such that 
$$
d\big(
X(t_i),\tilde \mM(c_i)
\big)
\leq \epsilon_i.
$$
If the cohomology classes $c_i$ are equal to a fixed
one $c_{\infty}$ for $i \geq i_0$,
then one can take $\eta$ such that 
$\eta |_ {[t'_{i_0},\infty)}$
is closed and periodic.
The trajectory  $X|_ {[t'_{i_0},\infty)}$
is then a minimizer of $L-c_{\infty}$, hence 
it is asymptotic to $\tilde \mL(c_{\infty})$
by definition.
The same holds for $\alpha$-limits.
\hfill

\newpage
\section{Globally minimizing orbits}\label{globally}
We have achieved our main goal, proving Theorem 
\ref{main}. 
However, 
the hypothesis of this theorem is rather abstract,
and some additional work is required in order 
to understand this hypothesis. In the present section,
we will describe the various sets of globally
minimizing orbits which have been defined in the literature.
Since most of the proofs have been written only in the 
autonomous case, we prove most of the results we state,
except the graph properties, 
mostly due to Mather,
and for which we send the reader
to \cite{Mather 1} and \cite{Mane}.

\subsection{}\label{Lcrit}
The Lagrangian $L$ is called critical if the infimum of the actions
of all periodic curves of all periods is $0$.
It is equivalent to require that the minimum of the actions
of all invariant probability measures is $0$.
Any Lagrangian satisfying the hypotheses of 
\ref{hypotheses}
can be made critical by the addition of a real constant.
See \ref{alpha} below for more details.
\subsection{}\label{functions}
Let $L$ be a critical Lagrangian.
For all $t'\geq t$ we define the function
\begin{align*}
F_{t,t'}: M\times M& \lto \Rm\\
 (x,x')  & \lmto \min _{\gamma\in \Gamma}
 \int _t ^{t'} L(d\gamma(u),u)\,du
\end{align*}
where the minimum is taken on the set $\Gamma$
of curves $\gamma\in C^1([t,t'],M)$ satisfying
$\gamma(t)=x$ and $\gamma(t')=x'$.
We also define, for each $(s,s')\in S^2$
the function
\begin{align*}
\Phi_{s,s'}: M\times M& \lto \Rm\\
 (x,x')  & \lmto \inf F_{t,t'}(x,x')
\end{align*}
where the infimum is taken 
on the set of $(t,t')\in \Rm^2$
such that $s=t\text{ mod }1$,  $s'=t'\text{ mod }1$, 
and $t'\geq t+1$. 
Following Mather, we introduce one more function 
\begin{align*}
h_{s,s'}: M\times M& \lto \Rm\\
 (x,x')  & \lmto \liminf_{t'-t\lto \infty}  F_{t,t'}(x,x')
\end{align*}
where the liminf
is restricted to the set of $(t,t')\in \Rm^2$
such that  $s=t\text{ mod }1$ and   $s'=t'\text{ mod }1$.
These 
functions have symmetric counterparts
$$
d_{s,s'}(x,x')=
h_{s,s'}(x,x')
+h_{s',s}(x',x)\;\;
\text{ and }\;\;
\tilde d_{s,s'}(x,x')=
\Phi_{s,s'}(x,x')
+\Phi_{s',s}(x',x)
$$
It is not hard to see, if $L$ is critical,
that  $d\geq \tilde d\geq 0$.

\subsection{}\label{F}
\textsc{Lemma :}
The set of function $F_{t,t'}$ with $t'\geq t+1$
is equilipschitz and equibounded.
\vs
\\
\textsc{Proof :}
Let us fix a number $\Delta \geq 1$ greater than the diameter of $M$.
In views of Lemma \ref{apriori}
and \ref{apriori2},
there exists a constant $K$ such that,
if $t'\geq t+1$ and if $\gamma\in C^1([t,t'],M)$
is a minimizer, then $\|d \gamma \|\leq K$.
Let us set 
$$
B=\max _{(z,t)\in TM\times \Rm, \|z\|\leq K+3\Delta}|L(z,t)|.
$$
Consider  $t'\geq t+1$ and four points
$x_0,x'_0,x_1,x'_1$ in $M$.
There is a minimizing  curve 
$\gamma\in  C^1([t,t'],M)$
such that
$A(\gamma)=F_{t,t'}(x_0,x'_0).$
Let us set 
$$
\delta=\min \{1/3,d(x_0,x_1)\}\;\;
\text{ and } \;\;
\delta'=\min \{1/3,d(x'_0,x'_1)\}.
$$
The geodesic $x\in C^1([t,t+\delta],M)$
between $x_1$ and $\gamma(t+\delta)$ satisfies
$$
\|dx\|\leq d(x_1,\gamma(t+\delta))/\delta
\leq \big( d(x_0,\gamma(t+\delta))+d(x_0,x_1) \big)/\delta
\leq
K+d(x_0,x_1)/\delta
\leq K+3\Delta,
$$
hence $A(x)\leq B\delta.$
The same estimate is true with the geodesic
$x'\in C^1([t'-\delta',t'],M)$ connecting
$\gamma(t'-\delta')$ and $x'_1.$ 
We have 
\begin{align*}
F_{t,t'}(x_1,x'_1)
&
\leq F_{t,t+\delta}(x_1,\gamma(t+\delta))
+
F_{t+\delta,t'-\delta'}(\gamma(t+\delta),\gamma(t'+\delta'))
+
F_{t'-\delta',t'}(\gamma(t'-\delta'),x'_1)
\\
&\leq
F_{t+\delta,t'-\delta'}(\gamma(t+\delta),\gamma(t'+\delta'))
+B\delta+B\delta'\\
&\leq 
F_{t,t'}(x_0,x'_0)
-
A(\gamma|_{[t,t+\delta]})
-A(\gamma|_{[t'-\delta,t']})
+B\delta+B\delta'
\\
&
\leq F_{t,t'}(x_0,x'_0)+2B\delta+2B\delta'.\\
&\leq  F_{t,t'}(x_0,x'_0) +2Bd(x_0,x_1)+2Bd(x'_0,x'_1).
\end{align*}
This proves that $2B$ is a Lipschitz constant for all
the functions $F_{t,t'}$ with $t'\geq t+1.$
We need to introduce some definition before we prove
that these functions are equibounded.
The proof will be given in \ref{equibound}.
\subsection{}
We have defined in \ref{tonelli}
two classes of orbits,
$\bar M$-minimizers and minimizers.
It is useful to define 
distinguished classes of minimizers.
Recall that $L$ is a critical Lagrangian.
A curve $\gamma\in C^1(I,M)$
is called
semi-static if 
$$
A\big(
\gamma|_{[a,b]}
\big)
= \Phi _{a\text{ mod 1},b\text{ mod 1}}
\big(\gamma(a),\gamma(b)
\big)
$$
for all $[a,b]\subset I$.
An orbit $X(t)=(d\gamma(t),t \text{ mod }1)$
is called semi-static if $\gamma$
is a semi-static curve.
It is clear that semi-static orbits are 
minimizing.
A curve $\gamma\in C^1(I,M)$
is called
static if 
$$
A\big(
\gamma|_{[a,b]}
\big)
= -\Phi _{b\text{ mod 1},a\text{ mod 1}}
\big(\gamma(b),\gamma(a)
\big)
$$
for all $[a,b]\subset I$.
If $\gamma$ is not semi-static, then there exists
$[a,b]$ such that 
$$
A\big(
\gamma|_{[a,b]}
\big)> \Phi _{a\text{ mod 1},b\text{ mod 1}}
\big(\gamma(a),\gamma(b)
\big)
$$
hence
$$
A\big(
\gamma|_{[a,b]}
\big)
+\Phi _{b\text{ mod 1},a\text{ mod 1}}
\big(\gamma(b),\gamma(a)
\big)>\tilde d_{a\text{ mod 1},b\text{ mod 1}}
\big(\gamma(a),\gamma(b)\big)\geq 0
$$
hence $\gamma$ is not static.
It follows that static curves are semi-static.
We call $\tilde \mN$ the union in $TM\times S$
of the images of global semi-static orbits
(semi-static orbits with $I=\Rm$)
and $\tilde \mA$ the union of global static orbits.
Clearly,
$$
\tilde \mA\subset \tilde \mN\subset \tilde \mB.
$$
It has became usual to call $\tilde \mA$ the Aubry set,
and $\tilde  \mN$ the Ma\~ne set. 

\subsection{}
\textsc{Lemma :}
We have the equivalence
$$
d_{s,s}(x,x)=0
\eq
\tilde d_{s,s}(x,x)=0
\eq
x\in \mA_s,
$$
and the set $\tilde \mA$ is a non empty compact invariant set.
\vs
\\
\textsc{Proof :}
Since $d\geq \tilde d\geq 0$, 
it is enough to prove that $d_{s,s}(x,x)=0$
if $\tilde d_{s,s}(x,x)=0$ to prove the first equivalence.
Assume that  $\tilde d_{s,s}(x,x)=0$.
Recall that 
 $
\tilde d_{s,s}(x,x)=
2 \Phi_{s,s}(x,x).
$
Either the infimum in the definition of $\Phi$
is a minimum, or it is a liminf.
If it is a liminf, the proof is over.
If it is reached, 
there is a curve $\gamma:[t,t']\lto M$
such that $\gamma(t)=\gamma(t')=x$
and $t\text{ mod 1}=s=t'\text{ mod 1},$
satisfying $A(\gamma)= 0.$
In this case, we can paste $\gamma$ with itself several times
and build a curve
$\gamma_k:[t,t+k(t'-t)]$
such that $\gamma_k(t)=\gamma_k(t+k(t'-t))=x$
and such that  $A(\gamma_k)= 0.$
It follows that 
$h_{s,s}(x,x)=0$,
hence  
$
d_{s,s}(x,x)=0.
$
This ends the proof of the first equivalence.

Let us suppose that $d_{s,s}(x,x)=0,$ and prove that 
$x\in \mA_s.$
There is a sequence $\gamma_k\in C^1([t_k,t'_k],M)$
of minimizing curves such that $A(\gamma_k)\lto 0$,
$\gamma_k(t_k)=x$, $\gamma_k(t'_k)=x$  and 
such that $t_k \text{ mod } 1=s=t'_k\text{ mod }1$
and $t'_k-t_k\lto \infty.$
By Lemma \ref{minlim}
we can suppose, taking a subsequence,
that the curves $x_k(t)=\gamma_k(t+[t_k])$
and $y_k(t)=\gamma_k(t+[t'_k])$
are converging uniformly on compact sets
to minimizers
$\gamma^+\in C^1([s,\infty),M)$
and $\gamma^-\in C^1((-\infty,s],M).$
In the above expressions, $[t]$ is the integer part of the 
real number $t$, and $s$ also denotes
the real number in $[0,1[$ such that $s\text{ mod } 1=s.$
Let $\gamma$ be the curve that coincides with
$\gamma^-$ and $\gamma^+$ on $(-\infty,s]$ and $[s,\infty)$.
Clearly, $\gamma(s)=x.$
If $t\leq s\leq s+1\leq t'$, then
\begin{eqnarray*}
&&A(\gamma|_{[t,t']}) +
\Phi_{t'\text{ mod }1,t\text{ mod }1}
(\gamma(t'),\gamma(t))\\
&=&
A(\gamma|_{[t,s]}) +A(\gamma|_{[s,t']})+
\Phi_{t'\text{ mod }1,t\text{ mod }1}
(\gamma(t'),\gamma(t))
\\
&=&
\lim 
\Big(
A\big(\gamma_k|_{[t+[t'_k],t'_k]}\big)
+A\big(\gamma_k|_{[t_k,[t_k]+t']}\big)\\
&+&\Phi_{t'\text{ mod }1,t\text{ mod }1}
(\gamma_k([t_k]+t'),\gamma_k([t'_k]+t))
\Big)\\
&\leq &\lim \Big(
A\big(\gamma_k|_{[t+[t'_k],t'_k]}\big)
+A\big(\gamma_k|_{[t_k,[t_k]+t']}\big)
+A\big(\gamma_k|_{[[t_k]+t',[t'_k]+t]}\big)
\Big)\\
&=&\lim A(\gamma_k)=0.
\end{eqnarray*}
hence the curve $\gamma$ is static and $x\in \mA_s.$
In order to
prove the last implication,
let us consider a static curve $\gamma.$
For each $t$, we have 
$$
\Phi_{t\text{ mod }1,t\text{ mod }1}
(\gamma(t),\gamma(t))
\leq
A\big( \gamma|_{[t-1,t+1]} \big)
+\Phi_{t\text{ mod }1,t\text{ mod }1}
(\gamma(t+1),\gamma(t-1))
=0.
$$
As a consequence, 
$$
\tilde d _{t\text{ mod }1,t\text{ mod }1}
(\gamma(t),\gamma(t))=0.
$$
Finally, the set $\mA$ is not empty
because it is clear that the minimum of the function
$x\lmto d_{s,s}(x,x)$ has to be $0$ for each $s$ if
$L$ is critical.
\hfill

\subsection{}\label{equibound}
We are now in a position to prove that
the functions $F_{t,t'}$, $t'\geq t+1$ are equibounded.
Let 
$$ 
A=\sup _{t, x, x'} F_{t,t+1/3}(x,x') 
\;\;\text{ and }\;\;
B=\sup_{s,s',x,x'} \Phi_{s,s'}(x,x'),
$$
both $A$ and $B$ are finite.
let $\gamma\in C^1(\Rm,M)$ be
a semi-static curve. There exist semi-static curves
since we just proved the existence of static curves.
Let us chose $t'\geq t+1$ and $x,x'\in M.$
We have
\begin{eqnarray*}
F_{t,t'}(x,x')  &\leq &
F_{t,t+1/3}\big(        x,\gamma(t+1/3)        \big)
+\\
&&
F_{t+1/3,t'-1/3}
\big(\gamma(t+1/3), \gamma(t'-1/3)
\big)
+F_{t'-1/3,t'}
\big(\gamma(t'-1/3), x'
\big)\\
& \leq &A+B+A,
\end{eqnarray*}
where we have used that, since  $\gamma$ is semi-static,
$$
F_{t+1/3,t'-1/3}
\big(\gamma(t+1/3), \gamma(t'-1/3)
\big)=
\Phi_{(t+1/3)\text{ mod }1,(t'-1/3)\text{ mod }1}
\big(
\gamma(t+1/3), \gamma(t'-1/3)
\big).
$$
Recalling that the functions $F_{t,t'}$ are equilipshitz,
we obtain the existence of a constant $C$ such that
$$
F_{t,t'}(x,x')\leq C
$$
for all $t'\geq t+1$ and all $(x,x')\in M^2$.
In order to end the proof, notice that,
when $k$ is large enough,
$$
F_{t,t'}(x,x')+F_{t',t+k}(x',x)\geq 0,
$$ 
hence
$F_{t,t'}\geq -C   .$
\hfill

\subsection{}
\textsc{Lemma :}
We have the inclusions
$$
\tilde \mM\subset
\tilde \mL\subset
\tilde \mA\subset
\tilde \mN\subset
\tilde \mB.
$$
\textsc{Proof :}
It is enough to prove that
$\tilde \mL \subset \tilde \mA.$
Let $X:[0,\infty)\lto TM\times S$ be a minimizing orbit
and let $\tilde \omega\in T_{\omega}M\times S$ be an omega-limit point.
Let $t_k\lto \infty $ be a sequence 
of times such that $X(t_k)\lto \tilde \omega$,
and assume that $s=t_k\text{ mod }1$ does not depend on $k$,
and that $t_{k+1}-t_k\lto \infty$.
Let $\gamma=\pi \circ X$.
Let us set $X_k(t)=X(t+[t_k])$.
Taking a subsequence if necessary,
we can suppose that the curves $X_k$ are
converging uniformly on compact
sets to a curve $Y(t)=(dx(t),t\text{ mod }1)$.
In order to prove that $x$ is a static curve,
we write, for $t'\geq t+1$,
\begin{eqnarray*}
&&A(x|_{[t,t']}) +
\Phi_{t'\text{ mod }1,t\text{ mod }1}
(x(t'),x(t))\\
&=&
\lim 
A\big(
\gamma|_{[t+[t_k],t'+[t_k]]}\big)
+\Phi_{t'\text{ mod }1,t\text{ mod }1}
(x(t'),x(t))\\
&=&
\lim 
\Big(
A\big(
\gamma|_{[t_{k-1},t_{k+1}]}
\big)-
A\big(
\gamma|_{[t_{k-1},t+[t_k]]}
\big)-
A\big(
\gamma|_{[t'+[t_k],t_{k+1}]}
\big)\Big)\\
&+&
\Phi_{t'\text{ mod }1,t\text{ mod }1}
(x(t'),x(t))\\
&\leq&
\liminf
\Big(
A\big(
\gamma|_{[t_{k-1},t_{k+1}]}
\big)\Big)
\\
&-&
\Big(
\Phi_{s,t\text{ mod }1}
(\omega,x(t))
+
\Phi_{t'\text{ mod }1,s}
(x(t'),\omega)
-
\Phi_{t'\text{ mod }1,t\text{ mod }1}
(x(t'),x(t))
\Big)\\
&\leq&
\liminf 
\Big(
A\big(
\gamma|_{[t_{k-1},t_{k+1}]}
\big)
\Big)
\leq 0.
\end{eqnarray*}
In this calculations,
we have used Lemma \ref{limit}
between the first line and the second,
and we have used Lemma \ref{F}
to obtain the last inequality.
More precisely,
it follows from this lemma
that the sum
$$
\sum_{k=1} ^n
A\big(
\gamma|_{[t_{2k-1},t_{2k+1}]}
\big)
=
A\big(
\gamma|_{[t_{1},t_{2n+1}]}
\big)
=
F_{t_{1},t_{2n+1}}\big(
\gamma(t_1),\gamma(t_{2n+1})
\big)
$$
is bounded, which implies that the liminf 
is not positive.
\hfill 
\subsection{}\label{fgp}
\textsc{First Graph property :}
Let us
call 
$\Pi:TM\times S\lto M\times S$
the natural projection.
Then $\Pi|_{\tilde \mA}$ is a bilipschitz
homeomorphism onto its image $\mA$.
In addition, we have
$$
\tilde \mN\cap \Pi^{-1}\big (\mA\big) =\tilde \mA.
$$
In other words, there is a Lipschitz section
$v:\mA\lto TM\times S$ of $\Pi$
with the property that,
for each $(x,s) \in \mA$, there is one and only one
semi-static orbit $X(t)$
satisfying $\Pi(X(0))=(x,s)$,
this orbit is static and is given by
$X(t)=\phi_t(v(x,s),s)$. 
%
%
%
%
%
%
%
%
%
\subsection{}\label{staticclasses}
It is not hard to see that 
$$
\tilde d_{s,s'}(x,x')=
d_{s,s'}(x,x')
$$
if $(x,s)\in \mA$ or $(x',s')\in \mA.$
We define an equivalence relation on $\mA$ by
saying that 
$(x,s)$ and $(x',s')$ are equivalent 
if and only if 
$d_{s,s'}(x,x')=0.$
We call static class an equivalence class of this relation.
We also call static class the image by 
the Lipschitz vector field $v$ 
of a static class in $M\times S.$
Static classes are compact invariant subsets of 
$\tilde \mA.$
\vs\\
\textsc{Remark :}
If $\gamma:[0,\infty)\lto M$
is minimizing, then the omega-limit
set of the orbit
$X(t)=(d\gamma,t\text{ mod }1)$
is contained in a static class.
\vs\\
\textsc{Proof :}
Let us consider sequences $t_k$ and $t'_k$ 
such that $t_k\text{ mod }1=s$ and $t'_k\text{ mod }1=s'$,
and such that $X(t_k)\lto \tilde \omega$
and $X(t'_k)\lto \tilde \omega'$.
We can assume in addition
that 
$t_k-t'_k\lto \infty$
and that $t'_k-t_{k-1} \lto \infty$.
If $\omega$ and $\omega'$ are the projections on
$M$ of $\tilde \omega$ and $\tilde \omega'$,
then
$$
d_{s,s'}(\omega,\omega')
\leq \liminf
A\big(
\gamma|_{[t_k,t_{k+1}]}
\big)
\leq \liminf 
F_{t_k,t_{k+1}}(\gamma(t_k),\gamma(t_{k+1}))
\leq 0,
$$
where the last liminf is not positive in view of Lemma \ref{F}
since $\gamma(t_k)$ is convergent.\hfill  \\
A semi-static curve thus 
has its alpha-limit contained in a static class,
and its omega-limit contained in a static class.
\vs\\
\textsc{Lemma :}
A semi-static curve is static if and only if
its alpha and omega-limit belong to the same 
static class. 
If $\tilde \mA$ contains only one static class,
then $\tilde \mN=\tilde \mA$.
It is the case for example if $\tilde \mM$ is 
transitive 
\textit{i.e.} if it has a dense orbit.
\vs \\
\textsc{Proof :}
It is quite clear that if $\gamma(t)$
is a static curve, then
$$
d_{t\text{ mod }1,t'\text{ mod }1}
(\gamma(t),\gamma(t'))=
\tilde d_{t\text{ mod }1,t'\text{ mod }1}
(\gamma(t),\gamma(t'))=0
$$
for all $t\leq t'$.
Taking the limit $t\lto -\infty$ and
$t'\lto \infty$
we obtain that the alpha and omega limit
belong to the same static class.
On the other hand, let $\gamma(t)$ be a semi-static curve
such that 
the alpha and omega-limit belong to the same static class.
Let us consider sequences
$t_k\lto -\infty$ and $t'_k\lto \infty$
of integers such that 
$\gamma(t_k)$ has a limit $\alpha\in M$ 
and $\gamma(t'_k)\lto \omega$.
The hypothesis is that $d_{0,0}(\alpha,\omega)=0$.
For each $t'\geq t$, we have
$$
d_{t\text{ mod }1,t'\text{ mod }1}
(\gamma(t),\gamma(t'))
+
d_{0,t\text{ mod }1}
(\alpha,\gamma(t))
+
d_{t'\text{ mod }1,0}
(\gamma(t'),\omega)
\leq d_{0,0}(\alpha,\omega)=0, 
$$
hence 
$d_{t\text{ mod }1,t'\text{ mod }1}
(\gamma(t),\gamma(t'))\leq 0
$
and $\gamma$ is static.\hfill  
\subsection{}\label{sgp}
If $\tilde \mS\subset  TM\times S$
is a static class,
we call $\tilde \mS^+$
the set of points $(z,s)\in TM\times \Rm$
such that the orbit $\phi_t(z,s)$
is semi-static on an open neighborhood 
of $[0,\infty)$, and omega-asymptotic to 
$\tilde \mS$.
We define $\tilde \mS^-$ in the same way with 
alpha-limits.
\vs
\\
\textsc{Second graph property :}
For each static class $\tilde \mS$,
the restriction of $\Pi$ to $\tilde \mS^+$
is a bilipschitz  homeomorphism onto its image, as well
as the restriction of $\Pi$ to $\tilde \mS^-$.
The set $\tilde \mN$ is the union of the graphs
$\tilde \mN\cap \tilde \mS^+$, as well as the union
of the graphs
$\tilde \mN\cap \tilde \mS^-$.

\subsection{}\label{alpha}
Let us now describe the action of adding 
a closed one-form to $L$ on the various sets we 
have defined.
We  identify $H^1(S,\Rm)$ with $\Rm$
in the standard way.
To a closed one-form $\eta$ on $M\times S$,
we associates the cohomology
$\lambda(\eta)$ of its restriction to 
$\{x\}\times S$, this cohomology does not depend
on $x\in M$, and depends only of the cohomology
of $\eta$.
Recall that we have defined in 
\ref{formes}
the class $[\eta]\in H^1(M,\Rm)$ of any closed
one form $\eta$ on $M\times S$.
the function
$$
\eta
\lmto ([\eta],\lambda(\eta))
$$
induces an isomorphism between 
$H^1(M\times S,\Rm)$ and $H^1(M,\Rm)\times\Rm$.
Let us fix a Lagrangian $L$, not necessarily critical.
We say that a closed one-form $\eta$ on $M\times S$
is critical if $L-\eta$ is critical.
\vs\\
\textsc{Theorem (Mather) :}
There exists a convex and superlinear function 
$$
\alpha:H^1(M,\Rm)\lto \Rm
$$
with the property that a closed  one-form
$\eta$ is critical
if and only if 
$$
\lambda(\eta)=-\alpha([\eta]).
$$
See \cite{Mather 1}  for the proof of this theorem
and for details on the following remarks.
The subderivative of $\alpha$ at a class $c$
is the set of rotation vectors in $H_1(M,\Rm)$
of the probability measures minimizing 
$L+c$.
It is usual to call
$$
\beta :H_1(M,\Rm)\lto \Rm
$$
the Fenchel transform of $\alpha$.
For each $\omega\in H_1(M,\Rm)$,
the number $\beta(\omega)$ is the minimal action of 
invariant probability measures of rotation vector $\omega$.
Given a critical form $\eta$, we can associate 
all the sets $\tilde \mM,\tilde\mA,\ldots$
to the critical Lagrangian $L-\eta$.
It is not hard to see
that these sets depend only on the class 
$[\eta]\in H^1(M,\Rm)$.
We define in the natural way the sets
$$
\tilde \mM(c)\subset
\tilde \mL(c)\subset
\tilde \mA(c)\subset
\tilde \mN(c)\subset
\tilde \mB(c)
$$
associated to the critical Lagrangian 
$L-\eta$, where $\eta$ is any critical form
satisfying $[\eta]=c$.
Notice that, in view of Mather's Theorem above,
the function $\eta\lmto [\eta]$ restricted
to critical forms is surjective.

\newpage
\section{Convergence of the Lax-Oleinik semigroup}\label{Lax}

The Graph properties provide a good description
of the Ma\~ne set $\tilde \mN$.
However, the set involved in the hypothesis of Theorem
\ref{main} is the
\textit{a priori } larger set $\tilde \mB$.
The relations between the sets  $\tilde \mB$ and 
 $\tilde \mN$ are related to the asymptotic
behavior of the so called Lax-Oleinik semi-group.
In all this section, we will consider a critical Lagrangian $L$
as  defined in \ref{Lcrit}. Results similar to the ones of this 
section have been obtained from the point of view of
Hamilton-Jacobi equations in \cite{Roquejoffre}.

\subsection{}\label{regular}
We say that $L$ is regular if the liminf in the definition
of the functions $h_{s,s'}$ 
given in \ref{functions} is a  limit for all
$s,s',x,x'$. 
In this case,
since the functions $F_{t,t'}$ are equilipschitz,
we have uniform convergence of 
the sequence $F_{t,t'}$, $t\text{ mod }1=s$, $t'\text{ mod }1=s'$  
to $h_{s,s'}$ for all $s,s'$.
If $L$ is regular and if $\eta$ is an exact one-form on $M\times S,$
then $L-\eta$ is regular.

\subsection{}\label{LO}
It is usual to define the mapping 
$T_t:C(M,\Rm)\lto C(M,\Rm)$ by
the expression
$$
T_t u(x)=\min _{y\in M}
\big(
u(y)+F_{0,t}(y,x)
\big).
$$
The sequence
$(T_n)_{n\in \Nm}$
is a semi-group called the Lax-Oleinik semi-group,
see \cite{Fathi} and \cite{FM}.
We say that 
the Lax-Oleinik semi-group is convergent
if, for each  function $u\in C(M,\Rm)$, there exists a function
$U\in C(M\times S,\Rm)$
such that 
$$
\lim _{t\text{ mod 1}=s, t\rightarrow \infty}
T_t u(x)
=
U(x,s).
$$
It is standard that the Lax-Oleinik semi-group is convergent
if and only if $L$ is regular, see \cite{Fathi} and \cite{FM}.
We shall recall the argument.
If $L$ is regular, then the Lax-Oleinik semi-group
is clearly convergent with limit
$$
U(x,s)=
\min_{y\in M}
\big(
u(y)+h_{0,s}(y,x)
\big).
$$
On the other hand, Assume that the Lax-Oleinik semi-group is
convergent.
Let us fix $t\in \Rm$ and $z\in M$,
and set $u(x)=F_{t,k}(z,x)$, where $k\in \Nm$ is chosen such that
$k\geq 1+t$.
For each $t'\geq k$, we have
$F_{t,t'}(z,x)=T_{t'-k}u(x)$.
If we fix $t'\text{ mod }1=s'$ and let $t'$ go to
infinity,
this is converging to $U(x,s')$,
which has to be equal to $h_{s,s'}(z,x)$.
It follows that $L$ is regular.

\subsection{}\label{NB}
\textsc{Proposition :}
If $L$ is regular,
then 
$
\tilde \mB=\tilde \mN.
$
\vs\\
\textsc{Proof :}
Let $\gamma \in C^1(\Rm,M)$ be a minimizing orbit.
We have to prove that this orbit is semi-static.
Let us consider a sequence  $t_k\lto -\infty$
such that $s=t_k\text{ mod }1$
does not depend on $k$ 
and such that 
$\alpha=\lim \gamma(t_k)$ exists.
In the same way, we consider a sequence $t'_k\lto \infty$
and set $s'=t'_k\text{ mod }1$
and $\omega=\lim \gamma(t'_k).$
We have 
$$
A\big(
\gamma|_{[t_k,t'_k]}
\big)
=
F_{t_k,t'_k}\big(
\gamma(t_k),\gamma(t'_k)
\big)
\lto 
h_{s,s'}
\big(\alpha,\omega\big).
$$
Let us consider a compact interval of times $[a,b]$,
and assume to make things
simpler that $s'= a\text{ mod }1$ and  $s= b\text{ mod }1$.
For $k$ large enough, we have
\begin{eqnarray*}
A\big(
\gamma|_{[a,b]}
\big)
&=&
A\big(
\gamma|_{[t_k,t'_k]}
\big)
-
A\big(
\gamma|_{[t_k,a]}
\big)
-
A\big(
\gamma|_{[b,t'_k]}
\big)\\
&=&
F_{t_k,t'_k}(\gamma(t_k),\gamma(t'_k))-
F_{t_k,a}(\gamma(t_k),\gamma(a))-
F_{b,t'_k}(\gamma(b),\gamma(t'_k)).
\end{eqnarray*}
Taking the limit,  we get 
$$
A\big(
\gamma|_{[a,b]}
\big)
=
h_{s,s'}(\alpha,\omega)
-
h_{s,s}(\alpha,\gamma(a))
-
h_{s',s'}(\gamma(b),\omega).
$$
On the other hand, we observe if $L$ is regular that
$$
h_{s,s'}(\alpha,\omega)
\leq 
h_{s,s}(\alpha,\gamma(a))
+
\Phi_{s,s'}(  \gamma(a),\gamma(b)   )+
h_{s',s'}(\gamma(b),\omega).
$$
As a consequence, we obtain
$$A\big(
\gamma|_{[a,b]}
\big)
\leq 
\Phi_{s,s'}(  \gamma(a),\gamma(b)   )
$$
hence $\gamma$ is semi-static.\hfill

\subsection{}\label{matherreg}
\textsc{Lemma :}
If for each $(x,s) \in \mM$,
the liminf in the definition of
$h_{s,s}(x,x)$
is a limit,
\textit{i.e.}
if
$$
F_{t,t+n}(x,x)\underset{n\rightarrow \infty}{\lto}0
$$
for each $(x,s)\in \mM$ and each $t$ satisfying
$t\text{ mod }1=s$,
then $L$ is regular.
\vs
\\
\textsc{Corollary :}
If $\tilde \mM$ is a union of 1-periodic orbits, then
$L$ is regular.
\vs
\\
\textsc{Proof :}
Let us fix $(x,s)$ and $(x',s')$ in $M\times S$, and $\epsilon>0$.
We want to prove that there exists $T$ such that,
if $t$ and $t'$ satisfy $t \text{ mod }1=s$, 
$t' \text{ mod }1=s'$ and $t'\geq t+T$,
then 
$$F_{t,t'}(x,x')\leq h_{s,s'}(x,x')+\epsilon.$$
Let $K$ be a common Lipschitz constant of all 
functions $F_{t,t'}$ with $t'\geq t+1$.
Such a constant exists by Lemma \ref{F}.
Let $\gamma:[t,t']\lto M$ be a minimizing curve
such that $A(\gamma)=F_{t,t'}(x,x')$
and such that $\gamma(t)=x$ and $\gamma(t')=x'$.
By Lemma \ref{sejour},
it is possible to chose 
$t_0\leq t_1\leq t'_0$
such that  $t_0 \text{ mod }1=s$ and 
$t'_0 \text{ mod }1=s'$, 
and a minimizing 
curve $ \gamma\in C^1([t_0,t'_0],M)$
such that $A(\gamma)=F_{t_0,t'_0}(x,x')$
and such that $\gamma(t_0)=x$,  $\gamma(t'_0)=x'$
and $d(\gamma(t_1),\mM_{t_1})\leq \epsilon/5K.$
Since $h_{s,s'}(x,x')=\liminf F_{t,t'}(x,x')$,
we can suppose in addition that 
$$
F_{t_0,t'_0}(x,x')\leq h_{s,s'}(x,x')+\epsilon/2.
$$
Let $x_1=\gamma(t_1)$, we have
$$
F_{t_0,t'_0}(x,x')
=
F_{t_0,t_1}(x,x_1)+
F_{t_1,t'_0}(x_1,x'),
$$
and there exists a point $y\in \mM_{t_1}$ such that
$d(x_1,y)\leq \epsilon/5K.$
It follows that 
$$
\big|
F_{t_0,t'_0}(x,x')
-
F_{t_0,t_1}(x,y)-
F_{t_1,t'_0}(y,x')
\big|
\leq\epsilon/2,
$$
hence
$$
F_{t_0,t_1}(x,y)+F_{t_1,t'_0}(y,x')
\leq h_{s,s'}(x,x')+\epsilon.
$$
Let us now consider 
$t$ and $t'$ such that  
$t \text{ mod }1=s$, 
$t' \text{ mod }1=s'$
and $t'-t=t'_0-t_0+n$ with $n\in\Nm$, we have 
$$
F_{t,t'}(x,x')=
F_{t_0,t'_0+n}(x,x')
\leq F_{t_0,t_1}(x,y)+
F_{t_1,t_1+n}(y,y)+
F_{t_1+n,t'_0+n}(y,x').
$$
Taking the limsup  yields
$$
\limsup
F_{t,t'}(x,x')
\leq F_{t_0,t_1}(x,y)
+0
+F_{t_1,t'_0}(y,x')
\leq h_{s,s'}(x,x')+\epsilon.
$$
Since this holds for all $\epsilon>0$, the lemma is proved.
Let us now prove the corollary.
If $\gamma\in C^1(\Rm,M)$ is 1-periodic
and minimizing,
then for each $t$ the sequence 
$$
F_{t,t+n}(\gamma(t),\gamma(t+n)
=nF_{t,t+1}(\gamma(t),\gamma(t+1))
$$
is bounded,
hence 
$
F_{t,t+n}(\gamma(t),\gamma(t))=0
$
for each $n$.
As a consequence, if $\tilde \mM$ is a union of
1-periodic orbits, then the hypothesis of the lemma is satisfied
and $L$ is regular.
\hfill

\subsection{}\label{rational}
One may wish to consider
the given Lagrangian $L$, which is 1-periodic in time,
as a $k$-periodic function of time only.
This is best done in our setting by
considering the new 1-periodic Lagrangian 
$$
L^k(x,v,t)=L(x,k^{-1}v,kt).
$$
This Lagrangian has the property that 
a curve $\gamma\in C^1(I,M)$ is an extremal
of $L^k$ if and only if
the curve $\gamma^k:t\lmto \gamma(kt)$
is an extremal of $L$.
We call $\mM^k$, $\mA^k$,\ldots
the various sets associated to $L^k$.
It is clear that 
$$\tilde \mB^k=\tilde \mB.$$
On the other hand, 
we have 
$$\tilde \mN \subset \tilde \mN^k$$
but 
it is not hard to build examples
where $\tilde \mN \neq \tilde \mN^k$
(see \cite{FM}).
Since $ \tilde \mN^k\subset\tilde \mB,$
this provides examples where 
$$
\tilde \mB\neq\tilde \mN.
$$
A direct consequence of Corollary \ref{matherreg}
and Proposition \ref{NB} is\\
\textsc{Lemma :}
If $\mM$ is a union of k-periodic orbits, then
$L^k$ is regular, hence
$
\tilde \mB=\tilde \mN^k.
$

\subsection{}\label{irrat}
\textsc{Lemma :}
If $\tilde \mM$ is minimal in the sense
of topological dynamics and
if there exists a sequence $\gamma_n$
of $n$-periodic curves such that 
$A(\gamma_n)\lto 0$, then 
$L$ is regular, hence
$\tilde \mA=\tilde \mN=\tilde \mB$.
\vs\\
\textsc{Proof :}
We can suppose that the curves $\gamma_n$ are minimizers.
Let us consider the n-periodic orbits
$X_n (t)=(d\gamma_n(t),t\text{ mod }1).$
Let us also note $X_n$ the image of $X_n$,
which is a compact subset of $TM\times S$.
Each subsequence of $X_n$ has a
convergent subsequence 
(for the Hau\ss dorff topology).
The limit of such a subsequence
is an invariant subset of $\tilde \mM$.
Since $\tilde \mM$ is minimal, this limit has to be 
$\tilde \mM$, hence $X_n$ is converging
to $\tilde \mM$ for the Hau\ss dorff topology.
It follows that each point
$(x,s)\in \mM$ is the limit of a sequence
$(\gamma_n(t_n),s)$ with $t_n\text{ mod }1 =s$
for each $n$.
Using Lemma \ref{F}, we get that
$$
\limsup F_{t,t+n}(x,x)
=\limsup F_{t,t+n}(\gamma_n(t_n),\gamma_n(t_n))
=\limsup A(\gamma_n)=0
$$
for each $(x,s)\in \mM$
and each $t$ satisfying $t\text{ mod }1=s$.
By Lemma \ref{matherreg}, $L$ is regular.
\hfill

\subsection{}
\textsc{Theorem (Fathi, \cite{Fathi}) :}
If $L$ does not depend on $t$, then it is regular.\\
As a consequence, in the autonomous case, the sets 
$\tilde \mB$ and $\tilde \mN$ are the same, hence our result is
precisely the result of Mather in this case. 
\newpage
\section{Twist Maps}\label{twist}
We are now going to specify our results in the case $M=S$.
As we shall see, unlike Mather's theorem of \cite{Mather},
our result in high dimension is optimal when restricted to this case,
in the sense that two cohomology classes 
$c$ and $c'$ are $C$-equivalent if and only if
the sets $\tilde \mB(c)$
and $\tilde \mB(c')$ are
not separated by a rotational invariant curve.

\subsection{}
Let $f:TS\lto TS$ be the 
Poincar\'e return map associated to the 
section $TS\times\{0\}$.
Moser has proved that any twist map
of the annulus $TS$  can be realized
as the Poincar\'e map of a Lagrangian flow
satisfying our hypotheses (\cite{Moser}).

\subsection{}
\textsc{Theorem :}
If $M=S$, then for each $c\in H^1(S,\Rm)$,
either $\tilde \mB(c)$
contains an invariant torus which is a Lipschitz graph,
or $R(c)=H^1(S,\Rm)$.

\subsection{}
We will now prove this theorem and give a
description of the invariant sets.
We identify $H^1(S,\Rm)$ and $H_1(S,\Rm)$ with $\Rm$
in the standard way. For each $c \in H^1(S,\Rm)$,
the set $\tilde \mA_0(c)$
is an $f$-invariant graph.
By the theory of homeomorphisms
of the circle,
the map $f$ restricted to  $\tilde \mA_0(c)$
has a rotation number, which is the only
subderivative of $\alpha$ at point $c$.
Hence $\alpha$ is differentiable, and $\alpha'(c)$
is the rotation number of $f|_{\tilde \mA_0(c)}$.

\subsection{}
\textsc{Irrational rotation number :}
Let us consider an irrational number $\omega$.
It is well known that $\beta$ is differentiable at 
$\omega$ (see \cite{Mather 3}) hence there exists only one value
$c$ such that $\alpha'(c)=\omega$.
It is clear that $\tilde \mM(c)$
is minimal in the sense of topological dynamics,
and we have
$$
\tilde\mA(c)=\tilde \mN(c)=\tilde \mB(c).
$$
As a consequence $\tilde \mB(c)$ is a graph.

\textsc{Proof :}
We can assume by adding a critical 
form $\eta$ satisfying $[\eta]=c$
that $\beta(\omega)=0=\beta'(\omega)$.
In view of Lemma $\ref{irrat}$,
it is enough to prove the existence of 
a sequence $\gamma_n$ of $n$-periodic orbits
such that $A_c(\gamma_n)\lto 0$.
Let us define $\gamma_n$ as the orbit 
minimizing the action $A(\xi|_{[0,n]})$
among the curves $\xi\in C^1(\Rm,S)$ 
whose liftings $\bar \xi$ to the universal cover $\Rm$
satisfy 
$\bar \xi(t+n)-\bar\xi(t)=[n\omega]$
for each $t$.
It is well known that 
$A(\gamma_n)=n\beta([n\omega]/n)$, which is converging to $0$
because  $\beta(\omega)=0=\beta'(\omega)$. \hfill

\subsection{}
\textsc{Rational rotation number :}
Let us consider a rational number $\omega=p/q$ in lowest terms.
Let us fix one $c$ such that $\alpha'(c)=\omega$.
The Mather set $\tilde \mM(c)$ 
is a union of $q$ periodic orbits.
By Lemma \ref{rational}, it follows that 
$\tilde \mB(c)= \tilde \mN^k(c)$.
Let $\mH$ be the closure of a connected component of the complement
of $ \mM(c)$ in $M\times S$.
The boundary of $\mH$ is made of two periodic curves
$\gamma^+$ and $\gamma^-$.
We see from the second graph property
that $\tilde \mB(c)\cap \Pi^{-1}\mH$
is the union of 
two graphs  $\tilde \mB^+$ and $\tilde \mB^-$,
where the orbits  $\tilde \mB^+$ are heteroclinic
from $\gamma^-$ to $\gamma^+$, as well as  
$\gamma^-$ and $\gamma^+$ themselves,
and the orbits of  $\tilde \mB^-$,
are heteroclinic
from $\gamma^+$ to $\gamma^-$ as well as $\gamma^-$ and $\gamma^+$.
If none of the projected sets
$\mB^+=\Pi(\tilde \mB^+)$ and $\mB^-=\Pi(\tilde \mB^-)$
is $\mH$, then their union is also  properly contained in $\mH$
\textit{i.e.}
$\mH\cap \mB(c)\neq \mH$.
In this case, $R(c)=\Rm$.
Else, $\tilde \mB(c)\cap \Pi^{-1}\mH$
contains a lipschitz graph.
If for all possible choice of $\mH$
the second option holds, then clearly
all the Lipschitz graphs can be glued together,
and $\tilde \mB(c)$ contains a Lipschitz graph.


\subsection{}
In terms of the Lax-Oleinik semi-group, we have 
proved the following.
Let $L$ be a critical Lagrangian,
and let $\omega$ be the rotation number of 
$\tilde \mA$.
Let us consider the integer $k$ 
defined by $k=1$ if $\omega$ is irrational,
and $k=q$ if $\omega=p/q$ in lowest terms.
Then the semi-group $T^k_n, n\in \Nm$ is converging.
Here we may view equivalently
$T^k_n$ as $T_{kn}$, 
or as the Lax-Oleinik semi-group
associated to $L^k$.
In other words,
the semi-group
$T_n$ has $k$-periodic asymptotic orbits.

\newpage
\section{The autonomous case}\label{autonomous}
We have seen that Theorem \ref{main}
is equivalent to the result stated by Mather in  \cite{Mather} 
in the autonomous case.
We shall now explain that this result is 
of no interest in the autonomous case.
I hope however that it is possible, 
still using the 
ideas introduced by Mather, to refine 
Theorem \ref{main}
in order to reach nontrivial applications
even in the autonomous case.

\subsection{}
A flat of $\alpha$ is
a closed convex  $K \subset H^1(M,\Rm)$
such that $\alpha|_K$ is an affine function.
To any closed convex set $K$, we associate the
vector subspace $VK=\text{Vect}( K-K)$.
A point $c$ is said to be in the interior of $K$ 
if there exists a neighborhood $U$ of $0$ 
in $VK$ such that $d+U\subset K.$
The interior of a flat is not empty.
Given $c\in  H^1(M,\Rm),$
we note $F(c)$ the union of all flats containing $c$
in their  interior.
It is clear that $F(c)$ is a flat,
we note  $VF(c)$ the associated vector space.

\subsection{}
Let $E(c)\subset  H^1(M,\Rm)$ be the vector subspace
of cohomology classes of one-forms of $M$ the support of
which are disjoint from $\mA(c)$.
Using the notations of \ref{R'}
we have 
$$
E(c)=V\big(   \mA(c)    \big)^{\bot}.
$$
In the autonomous case, we clearly have
$$
R(c)=V\big(   \mB(c)    \big)^{\bot}\subset E(c)
$$
since $\mA(c)\subset \mB(c).$
On the other hand, Massart \cite{Massart} has proved 
that $E(c)\subset VF(c)$, hence
$$
R(c)\subset E(c)\subset VF(c).
$$
From this follows that any admissible curve $c(t)$ is contained 
in a flat of $\alpha.$
Hence each C-equivalence class is contained in a flat.

\subsection{}
If $F$ is a flat of $\alpha$, there exists an Aubry set
$\mA(F)$ which is the aubry set $\mA(c)$ for 
all cohomology class $c$ in the interior of $F$,
and is contained in the aubry set of any cohomology class
$c\in F$.
This is also proved in \cite{Massart}.
As a consequence, there exists a Mather set $\mM(F)$
which is contained in all the Mather sets of the cohomology classes
of the flat.

\subsection{}
Let $C$ be a C-equivalence class.
It is contained in a maximal flat $F$.
It is not hard to see that the orbit
$(d\gamma(t),t \text{ mod }1)$ satisfies all 
the conclusions of Theorem \ref{main} 
if $\gamma(t) \in \mM(F).$
It follows that Theorem  \ref{main}
is of no interest in the autonomous case.

\end{document}